\documentclass[11pt,leqno]{article} 
\usepackage{graphics}
\newtheorem{thm}{Theorem}[section]

\newtheorem{prop}{Proposition}

\newcommand{\beqa}{\begin{eqnarray}}
\newcommand{\eeqa}{\end{eqnarray}}

\newcommand{\pf}{\noindent {\bf Proof:} $\s$ }
\newcommand{\epf}{ \hfill$\diamondsuit$ \medskip}

\newcommand{\beq}{\begin{equation}}
\newcommand{\eeq}{\end{equation}}
\newcommand{\lbl}{\label}
\newcommand{\s}{\; \;}

\newcommand{\ra}{\rightarrow}

\title{On separated solutions of logistic population equation with harvesting}

\author{
Philip Korman   \\ 
Department of Mathematical Sciences \\ 
University of Cincinnati \\ 
Cincinnati Ohio 45221-0025 \\
}

\date{}

\begin{document}

\maketitle
\begin{abstract} 
We provide a surprising answer to a question raised in  S. Ahmad and A.C. Lazer \cite{AS}, and extend the results of that paper.
 \end{abstract}

\begin{flushleft}
\end{flushleft}

\begin{flushleft}
\end{flushleft}

\section{An example}
\setcounter{equation}{0}
\setcounter{thm}{0}
\setcounter{lma}{0}

If one considers the  logistic population model 
\beq
\lbl{1}
z'(t)=a(t)z(t)-z^2(t) \,,
\eeq
with $a(t)$ continuous, positive and periodic function, of period $T$, then by a straightforward integration of this Bernoulli's equation one shows that there exists a unique positive solution $z_0(t)$ of period $T$, which attracts all other positive solutions as $t \ra \infty$, see e.g., M.N. Nkashama \cite{N}. The problem also has the trivial solution $z=0$. To account for  harvesting, one may consider the model 
\beq
\lbl{2}
z'(t)=a(t)z(t)-z^2(t)-k\gamma(t) \,,
\eeq
where $\gamma(t)$ is also  a continuous, positive and periodic function, of period $T$, and $k>0$ is a parameter. It was shown in \cite{A} and \cite{M} that there exists a $\bar k$, so that for $0<k<\bar k$ the problem (\ref{2}) has exactly two positive solutions  of period $T$, exactly one positive $T$-periodic solution at $k=\bar k$, and no $T$-periodic solutions for $k>\bar k$. Moreover, one can show that there is a curve of solutions beginning with  $z_0(t)$ at $k=0$ which is decreasing in $k$, and a curve of solutions beginning with  $z=0$ at $k=0$ which is increasing in $k$. At $k=\bar k$ these solutions coincide, and then disappear for  $k>\bar k$.
\medskip

In a very interesting recent paper S. Ahmad and A.C. Lazer \cite{AS} studied the equation (\ref{2}) without the periodicity assumption on $a(t)$ and $\gamma (t)$. They introduced the key concept of {\em separated solutions} to take place of the periodic ones. Two  solutions $z_1(t)>z_2(t)$ are called  separated if
\[
\int _0^{\infty} \left[z_1(t)-z_2(t) \right] \, dt=\infty \,.
\]
They proved the following result (among a number of other results).

\begin{thm}\lbl{thm:1} (\cite{AS})
Let the functions  $a(t)$ and $\gamma (t)$ be continuous and bounded by positive constants from above and from below on $[0,\infty)$. Then there exists a critical $\bar k$, so that for $0<k<\bar k$ the problem (\ref{2}) has  two separated  positive solutions. At $k=\bar k$ there exists a positive bounded solution,  while for $k>\bar k$ there are no bounded positive solutions.
\end{thm}

The authors of \cite{AS} asked if it is possible for two separated solutions to exist at  $k=\bar k$. Our next example shows that the answer is yes, which appears to be counter-intuitive.
\medskip

We consider the equation 
\[
z'(t)=z(t)-z^2(t)-\frac14 k \,,
\]
i.e., $a(t)=1$, $\gamma(t)=\frac14$. When $0<k<1$ this equations has two constant solutions, which are the roots of the quadratic equation $z-z^2-\frac14 k=0$. (By \cite{AS} there are no  solutions separated from both of the constant ones.) At $k=1$ there are two separated solutions $z=\frac12$ and $z=\frac12+\frac{1}{t+1}$, while for $k>1$ all solutions go to $-\infty$ in finite time, as can be seen by writing this equation in the form
\[
z'(t)=-\frac14 (2z-1)^2-\frac14 (k-1) \,.
\]
Here $\bar k=1$.

\section{Separated from zero solution of the logistic equation}
\setcounter{equation}{0}
\setcounter{thm}{0}
\setcounter{lma}{0}

We consider now the logistic model ($ t\geq 0$)
\beq
\lbl{4}
z'(t)=a(t)z(t)-b(t)z^2(t) \,, \s z(0)>0 \,.
\eeq
There is a zero solution $z=0$. Any solution of (\ref{4}) is positive, and by definition it is   separated from zero if $\int_0^{\infty} z(t) \, dt=\infty$.

\begin{prop}
\label{prop:1}
Let the functions  $a(t)$ and $b(t)$ be continuous and satisfy $a(t) \leq A$, $0<b_1<b(t)<b_2$ on $[0,\infty)$, with positive constants $A$, $b_1$, $b_2$. Then any solution of (\ref{4}) is bounded and separated from zero if and only if
\beq
\lbl{10}
\int_0^{\infty} b(t) e^{\int  _0^t a(s) \, ds} \, dt=\infty \,.
\eeq
\end{prop}

\pf
Setting $1/z=v$ and $\mu (t)=e^{\int  _0^t a(s) \, ds} $, we integrate (\ref{4}) to obtain (here $c=\frac{1}{z(0)}$)
\beq
\lbl{11}
z(t)=\frac{\mu(t)}{c+\int  _0^t b(s) \mu (s) \, ds}=\frac{1}{b(t)} \frac{g'(t)}{g(t)} \,,
\eeq
with $g(t)=c+\int  _0^t b(s) \mu (s) \, ds$. Then
\[
\frac{1}{b_2} \left(\ln g(\infty)-\ln g(0) \right) \leq  \int_0^{\infty} z(t) \, dt  \leq  \frac{1}{b_1} \left(\ln g(\infty)-\ln g(0) \right) \,,
\]
and the proof follows.
\epf

\begin{prop}
\lbl{prop:10}
Let the functions  $a(t)$ and $b(t)$ be continuous and satisfy $a(t) \leq A$, $0<b_1<b(t)<b_2$ on $[0,\infty)$, with positive constants $A$, $b_1$, $b_2$. Assume that $J \equiv \int_0^{\infty} b(t) e^{\int  _0^t a(s) \, ds} \, dt<\infty$, and $a(t) \leq 0$ for large $t$. Then all  solutions of (\ref{4}) tend to zero as $t \ra \infty$.
\end{prop}

\pf
We have $\mu '=a(t)\mu \leq 0$ for large $t$. Since $J<\infty$, it follows that $\mu (t) \ra 0$ as $t \ra \infty$. Then $z(t) \ra 0$ by (\ref{11}).
\epf

The situation is different in case of negative initial data:
\beq
\lbl{12}
z'(t)=a(t)z(t)-b(t)z^2(t) \,, \s z(0)<0 \,.
\eeq

\begin{prop}
\label{prop:2}
Case $1$. 
If the condition (\ref{10}) holds, then all solutions of (\ref{12}) go to $-\infty$ in finite time. 
\medskip

\noindent
Case $2$. If
\[
J \equiv \int_0^{\infty} b(t) e^{\int  _0^t a(s) \, ds} \, dt<\infty \,,
\]
then solutions with $z(0)<-\frac1J$ go to $-\infty$ in finite time, while solutions with $z(0) \in (0,-\frac1J]$ exist for all $t>0$, and the solution with $z(0)=-\frac1J$ is separated from zero. Moreover, under the additional assumptions that $a(t) \leq 0$ for large $t$ and $\lim _{t \ra \infty} a(t)=0$, solutions with $z(0) \in (0,-\frac1J]$ tend to zero as $t \ra \infty$.
\end{prop}

\pf
Solutions of (\ref{12}) are given by (\ref{11}), with $c=\frac{1}{z(0)}<0$, from which the claims on blow up and global existence follow. When  $z(0)=-\frac1J$, we have
\beq
\lbl{25}
z(t)=\frac{\mu (t)}{-J+\int_0^t b(s) \mu (s) \, ds} \,,
\eeq
and then 
\[
\int _0^{\infty} z(t) \,dt \leq \frac{1}{b_2} \int _0^{\infty}  \frac{b(t)\mu (t)}{-J+\int_0^t b(s) \mu (s) \, ds} \, dt
\]
\[
=\frac{1}{b_2} \ln |-J+\int_0^t b(s) \mu (s) \, ds| {\large |}_{_0}^{^{\infty}}=-\infty \,,
\]
so that $\int _0^{\infty} z(t) \,dt=-\infty$, and $z(t)$ is separated from zero. The proof that $z(t) \ra 0$ follows as in Proposition \ref{prop:10} in case $z(0) \in (0,-\frac1J)$, and by L'Hospital's rule for $z(0)=-\frac1J$.
\epf

If $p(t)$ is any particular solution of (\ref{2}), defined for $t \in [0,\infty)$, and $z(t)$ is any other solution of  (\ref{2}), then $v(t)=z(t)-p(t)$ satisfies Bernoulli's equation
\beq
\lbl{9}
v'=\left[a(t)-2p(t) \right] v-v^2 \,.
\eeq
It follows that  any other solution of  (\ref{2}), which is larger than $p(t)$, is separated from $p(t)$ if and only if
\beq
\lbl{5}
\int_0^{\infty}  e^{\int  _0^t \left[a(s)-2p(s) \right] \, ds} \, dt=\infty \,.
\eeq

The anonymous reviewer of this paper posed the following question.  
\medskip

\noindent
{\bf Question}. Are there always separated solutions at $\bar k \,$?
\medskip

We shall show  that the answer is affirmative. We consider first  the periodic case, where the picture  is simpler.
\begin{prop}
In the conditions of the Theorem \ref{thm:1} (from \cite{AS}),  assume additionally that $a(t)$ and $\gamma(t)$ are $T$-periodic, i.e., $a(t+T)=a(t)$ and $\gamma(t+T)=\gamma(t)$ for all $t \in [0,\infty)$, and some $T>0$. Let $p(t)$ be the unique  $T$-periodic  solution of  (\ref{2})  at $k=\bar k$. Then   any other solution of  (\ref{2}), which is larger than $p(t)$, is bounded on $[0,\infty)$, and is separated from $p(t)$.
\end{prop}

\pf
As we mentioned above, there exists a $\bar k$, so that for $0<k<\bar k$ the equation  (\ref{2}) has exactly two positive solutions  of period $T$, exactly one positive $T$-periodic solution at $k=\bar k$, and no $T$-periodic solutions for $k>\bar k$. So that $(\bar k,p(t))$ is a ``turning point"  of $T$-periodic  solutions of  (\ref{2}). It follows that the corresponding linearized problem
\beq
\lbl{87}
w'(t)=\left[a(t)-2p(t) \right] w(t) \,, \s\s w(t+T)=w(t)
\eeq
has non-trivial solutions, which happens if and only if
\beq
\lbl{6}
\int_0^T \left[a(t)-2p(t) \right] \, dt=0 \,.
\eeq
(If (\ref{87}) had only the trivial solution, the  equation  (\ref{2}) would have $T$-periodic  solution for $k> \bar k$, by the implicit function theorem, see e.g., \cite{K} for similar arguments.) We claim that there is a constant $\alpha$, such that 
\[
\int  _0^t \left[a(t)-2p(t) \right] \, dt \geq \alpha \,, \s \mbox{for all $t>0$} \,.
\]
Indeed, for any $t>0$, we can find an integer $n \geq 0$, so that $nT \leq t \leq (n+1)T$. Using (\ref{6}),
\[
\int  _0^t \left[a(t)-2p(t) \right] \, dt=\int  _{nT}^t \left[a(t)-2p(t) \right] \, dt \,,
\]
and the integral of the periodic function $a(t)-2p(t)$ over an interval of length $<T$, is bounded below by some constant $\alpha$. Then (\ref{5}) holds, and the proof follows.
\epf

We now consider the general case.

\begin{thm}
Let $p(t)$ be a positive bounded solution of (\ref{2}) at $k=\bar k$ from the Theorem \ref{thm:1}. Calculate $I=\int _0^{\infty} e^{\int_0^t \left(a(s)-2p(s) \right) \, ds} \, dt$.
\medskip

\noindent
Case $1$. $I=\infty$. Then all solutions  of (\ref{2}), lying above $p(t)$ (i.e., $z(0)>p(0)$)  are bounded for all $t>0$ and separated from  $p(t)$, while all solutions below  $p(t)$ go to $-\infty$ in finite time.
\medskip

\noindent
Case 2. $I<\infty$, $a(t)-2p(t) \leq 0$ for large $t$ and $\lim _{t \ra \infty} \left(a(t)-2p(t) \right) =0$. Then for $z(0) \in [p(0)-\frac1I,\infty)$ solutions tend to $p(t)$, and at $z(0) = p(0)-\frac1I$ one has a positive bounded solution, call it $p_1(t)$, which is separated from $p(t)$ (and tending to $p(t)$).
\end{thm}

\pf
Since $v(t)=z(t)-p(t)$ satisfies (\ref{9}), all of the claims, except for positivity of $p_1(t)$ follow by the Propositions \ref{prop:1}, \ref{prop:10} and \ref{prop:2}.
If $p_1(t)$ was negative at some point, we could find a point $t_0$ at which $p_1(t_0)=0$ and $p'_1(t_0) \geq 0$, since $p_1(t)$ tends to $p(t)>0$. But then we have a contradiction in the equation   (\ref{2})  at $t=t_0$.
\epf

Both of the cases described in this theorem actually occur. The periodic equations, considered above, provide  examples for the first case. The second case occurs in the following example.
\medskip

\noindent
{\bf Example} Consider
\beq
\lbl{20}
z'=z-z^2-k \left(\frac14-\frac{2}{(t+5)^2} \right) \,,
\eeq
with $f(t)=\frac14-\frac{2}{(t+5)^2}>0$. At $k=1$ there is a bounded solution $p(t)=\frac12+\frac{2}{t+5}$. To see that $\bar k=1$, we write this equation as 
\[
z'=-\left(z-\frac12 \right)^2-\frac{k-1}{4}+\frac{2k}{\left(t+5 \right)^2} \,.
\]
Then for $k>1$, the right hand side is smaller than, say $-\frac{k-1}{8}$ for $t$ large, and hence all solutions go to $-\infty$ in finite time. Here $a(t)=1$, $a(t)-2p(t)=-\frac{4}{t+5}$. Compute $\mu (t)=e^{-\int _0^t \frac{4}{s+5} \, ds} =\frac{5^4}{(t+5)^4}$, $\int _0^{\infty} \mu (t) \, dt=\frac53$.
We have $p_1(t)=p(t)+v(t)$, where by (\ref{25})
\[
v(t)=\frac{\mu (t)}{-\int _0^{\infty} \mu (s) \, ds+\int _0^t \mu (s) \, ds}=\frac{\mu (t)}{-\int _t^{\infty} \mu (s) \, ds}=-\frac{3}{t+5} \,.
\]
Hence, at $k=1$, $p_1(t)=\frac12-\frac{1}{t+5}$ is the bounded positive solution, separated from $p(t)$.  At $k=1$, any solution of (\ref{20}), with $z(0) \in [\frac{3}{10}, \infty)$ tends to $p(t)$ as $t \ra \infty$, while solutions with $z(0)<\frac{3}{10}$ go to $-\infty$ in finite time ($p_1(0)=\frac{3}{10}$).
\medskip

\noindent
{\bf Acknowledgment:} I wish to thank the referee for posing a stimulating question.


\begin{thebibliography}{99}
\bibitem{A}
S. Ahmad, Convergence and ultimate bounds of solutions of nonautonomous Volterra-Lotka competition equations, {\em J. Math. Anal. Appl.} {\bf 127 (2)},  377-387 (1987).
\bibitem{AS}
S. Ahmad and A.C. Lazer, Separated solutions of logistic equation with nonperiodic harvesting, {\em J. Math. Anal. Appl.} {\bf  445}, no. 1, 710-718  (2017).

\bibitem{K}
P. Korman, Global Solution Curves for Semilinear Elliptic Equations, World Scientific, Hackensack, NJ (2012).


\bibitem{M}
J. Mawhin, First order ordinary differential equations with several periodic solutions, {\em J. Appl. Math. Phys.} {\bf  38}, 257-265   (1987). 

\bibitem{N}
M.N. Nkashama, Dynamics of logistic equations with non-autonomous bounded coefficients, {\em  Electron. J. Differential Equations}, No. 2 (2000).

\end{thebibliography}
\end{document}